\documentclass[11pt,a4paper]{article}

\usepackage[T1]{fontenc}
\usepackage[utf8]{inputenc}
\usepackage{authblk}

\usepackage{latexsym}
\usepackage{amsmath}
\usepackage{graphicx}

\usepackage{amssymb}
\usepackage{epsfig}

\usepackage{amsfonts} 
\usepackage{url}      

\def\be{\begin{equation}}
\def\ee{\end{equation}}

 \newenvironment{keywords}{ \noindent {\small \textbf{Key Words}}.}{ }

\begin{document}


\newcommand{\nms}{\normalsize}

\title{Deterministic approaches for solving practical black-box global optimization
problems}

\author[1,2]{Dmitri E. Kvasov\thanks{email: kvadim@dimes.unical.it }\thanks{Corresponding author. }}
\author[1,2]{ Ya.D. Sergeyev\thanks{ email: yaro@dimes.unical.it }
}
 \nms
 \affil[1]{\nms DIMES, Universit\`a della  Calabria, Italy}
\affil[2]{\nms Department of Software and Supercomputing
Technologies\\ Lobachevsky State University of Nizhni Novgorod,
Russia}

\date{}

\maketitle

\begin{abstract}
In many important design problems, some decisions should be made by
finding the global optimum of a multiextremal objective function
subject to a set of constrains. Frequently, especially in
engineering applications, the functions involved in optimization
process are black-box with unknown analytical representations and
hard to evaluate. Such computationally challenging decision-making
problems often cannot be solved by traditional optimization
techniques based on strong suppositions about the problem
(convexity, differentiability, etc.). Nature and evolutionary
inspired metaheuristics are also not always successful in finding
global solutions to these problems due to their multiextremal
character. In this paper, some innovative and powerful
deterministic approaches developed by the authors to construct
numerical methods for solving the mentioned problems are surveyed.
Their efficiency is shown on solving both the classes of random
test problems and some practical engineering tasks.
\end{abstract}




\section{Introduction}

Numerical approaches to efficiently find optimal parameters of
mathematical models arising from different real-life design
problems are nowadays becoming more and more significant in
industrial processes. Optimization models characterized by the
functions with several local optima (typically, their number is
unknown and can be very high) have a particular importance for
practical applications. When the best set of parameters should be
determined for these multiextremal models, traditional local
optimization techniques, including many heuristic approaches, can
be insufficient and, therefore, global optimization methods are
used. Moreover, the objective functions and constraints to be
examined are often black-box and hard to evaluate functions with
unknown analytical representations. For example, their values can
be obtained by executing some computationally expensive
simulation, by performing a set of experiments, and so on. Such a
kind of functions is frequently met in various fields of human
activity (as, e.g., automatics and robotics, structural
optimization, safety verification problems, engineering design,
network and transportation problems, mechanical design, chemistry
and molecular biology, economics and finance, data classification,
etc.) and corresponds to computationally challenging global
optimization problems, being actively studied around the world
(see, e.g., \cite{Audet:et:al.(2005), Floudas&Pardalos(2009),
Grossmann(1996), Horst&Pardalos(1995), Michalewicz(1996),
Mockus(2000), Pardalos&Romeijn(2002),
Paulavicius&Zilinskas(2013b), Pinter(1996), Sergeyev&Kvasov(2008),
Sergeyev&Kvasov(2011), Strongin&Sergeyev(2000),
Zhigljavsky&Zilinskas(2008)} and the references given therein).

This paper is based upon the work \cite{Kvasov&Sergeyev(2012c)}
presented at the Eleventh International Conference on
Computational Structures Technology (Dubrovnik, Croatia, 4--7
September 2012) and extends the previous research in both the
theoretical and the experimental directions, thus offering to the
optimal design community competitive tools to tackle real-life
engineering decision-making tasks.

The paper is structured as follows. In the next Section, an
insight into black-box global optimization problems is given and
some approaches for their solution are briefly discussed. One of
these approaches are based on a quite natural (from the physical
viewpoint) assumption that the objective function and constraints
have bounded slopes, i.e., they satisfy the Lipschitz condition.
Such a problem statement is formalized and examined more in detail
in Section~3. Some approaches proposed by the authors to construct
efficient numerical methods for solving the mentioned problems are
briefly presented in Section~4. Section~5 illustrates the
theoretical considerations of the paper with some numerical
experiments. Finally, conclusions and future research directions
are drawn.

\section{Black-box global optimization}

To illustrate real-world global optimization problems we deal
with, let us consider the field of geomechanics and geophysics
where one has to work with different mechanical-mathematical
optimization models. Generally, these models are very complex
since they can involve multidimensional linear or nonlinear
partial differential equations with multiple contact boundaries,
regions with sharp changes in functions values, ill-posedness, and
so on. Knowledge of the properties and types of geological rocks
lying at a depth of several kilometers is of great interest, e.g.,
for prospecting seismology, which determines the location of oil
fields by means of acoustic waves. Prospecting seismology
techniques allow one both to avoid the costly exploration methods
(e.g., well drilling) and to accelerate the process of pinpointing
oil resources. Among these techniques, numerical methods for
solving inverse problems are of fundamental importance in
prospecting seismology. They aim at estimating parameters of the
Earth's structure and material properties (e.g., location of the
inhomogeneities as cracks/cavities in the crust) based on data
measured on the surface.

To give a concrete example of the resulting global optimization
problem, a simplified version of the prospecting seismology
inverse problem can be taken: namely, let us suppose that there is
a fluid-filled crack (or a number of such cracks) of a given
length located in the host rock with known elastic properties
(see, e.g., \cite{Golubev:et:al.(2010), Kvasov:et:al.(2012)}).
Then, the vector $x$ of unknown parameters defining the region
geometry contains only two components: the depth of the crack
occurrence $h$, $h_1 \leq h \leq h_2$, and the crack inclination
angle $\alpha$, $\alpha_1 \leq \alpha \leq \alpha_2$ ($h_1$,
$h_2$, $\alpha_1$, and $\alpha_2$ are known constants).

One of the peculiarities of the stated problem is that information
can be obtained only from experimental measurements with the usage
of acoustic sounding (see, e.g., \cite{Leviant:et:al.(2007)}). A
number of seismic detectors are located at points $d_i$ on the
surface of the Earth, which record the vertical components
$\tilde{V_y}(d_i, t_j)$ of particles velocity in the reflected
wave at time instances $t_j$. We look for such a value of $x$ that
best fits the numerically simulated response ${V_y}(x, d_i, t_j)$
to a measured one. Computational simulation can be performed by
some numerical integration algorithm: for example, the
grid-characteristic method (see, e.g.,
\cite{Petrov&Kholodov(1984), KvasovI&Petrov(2012)}) can be used
for this scope, thus taking into account the physical features of
the problem and allowing one to set correctly the boundary and
contact conditions.

Hereby, this particular problem can be formulated as the following
least squares optimization problem (see, e.g.,
\cite{Strongin&Sergeyev(2000), Zhigljavsky&Zilinskas(2008),
Ljung(1999), Schittkowski(2002), Sen&Stoffa(1995),
Stavroulakis(2001)}):

\begin{equation}
 \min f(x), \hspace{3mm} x \in D = [h_1,h_2] \times [\alpha_1,
 \alpha_2], \label{eq_app1}
\end{equation}

\begin{equation}
 f(x) = \sum_i \sum_j [ V_y (x, d_i, t_j) -
 \tilde{V_y}(d_i,t_j)]^2 \label{eq_app2}.
\end{equation}

Function~\eqref{eq_app2} is essentially multiextremal, it has no
analytical representation and its evaluation (sometimes called
trial) is associated with performing computationally expensive
numerical experiments. Therefore, the usage of fast and robust
global optimization methods aimed at tackling this class of
complex multiextremal problems is required for solving efficiently problem~\eqref{eq_app1}--\eqref{eq_app2}.

Because of the computational costs involved, typically a small
number of functions evaluations are available for a decision-maker
(engineer, physicist, chemist, economist, etc.) when optimizing
such costly functions. Thus, the main goal is to develop fast
global optimization methods that produce acceptable solutions with
a limited number of functions evaluations. But to obtain this
goal, there are still a lot of difficulties that are mainly
related either to the lack of information about the objective
function (and constraints, if any) or to the impossibility to
adequately represent the available information about the
functions.

For example, gradient-based algorithms (see, e.g.,
\cite{Floudas&Pardalos(2009), Horst&Pardalos(1995),
Pardalos&Romeijn(2002)}) cannot be used in many applications
because black-box functions are often non-differentiable or
derivatives are not available and their finite-difference
approximations are too expensive to obtain. Automatic
differentiation (see, e.g., \cite{Corliss:et:al.(2002)}), as well
as interval-based approaches (see, e.g., \cite{Csendes(2000),
Kearfott(1996)}), cannot be appropriately used in cases of
black-box functions when their source codes are not available. A
simple alternative could be the usage of the so-called direct (or
derivative-free) search methods (see, e.g., \cite{Kelley(1999),
Conn:et:al.(2009), Kolda:et:al.(2003), Liuzzi:et:al.(2010),
Custodio:et:al.(2011), Rios&Sahinidis(2013)}), frequently used now
for solving engineering design problems (see, e.g., the DIRECT
method \cite{Floudas&Pardalos(2009), Kelley(1999),
Jones:et:al.(1993)}, the response surface, or surrogate model
methods~\cite{Jones:et:al.(1998), Regis&Shoemaker(2005)}, pattern
search methods~\cite{Torczon(1997), Custodio:et:al.(2007),
Audet&Dennis(2003), Audet&Dennis(2006)}, etc.). But unfortunately
(see, e.g.,~\cite{Sergeyev&Kvasov(2006), DiSerafino:et:al.(2011),
Kvasov&Sergeyev(2012a), Paulavicius&Zilinskas(2013a)}), these
methods either are designed to find only stationary points or can
require too high computational effort for their work.

Therefore, solving the described global optimization problems is
actually a real challenge both for the theoretical and the applied
scientists. In this context, deterministic global optimization is
a well developed mathematical theory which has many important
applications (see, e.g., \cite{Floudas&Pardalos(2009),
Horst&Pardalos(1995), Pinter(1996), Sergeyev&Kvasov(2008),
Strongin&Sergeyev(2000), Floudas(2000)}). One of its main
advantages is the possibility to obtain guaranteed estimations of
global solutions and to demonstrate (under certain analytical
conditions) rigorous global convergence properties. However, the
currently available deterministic models can still require too
large number of functions evaluations to obtain adequately good
solutions for these problems.

Stochastic approaches (see, e.g., \cite{Floudas&Pardalos(2009),
Horst&Pardalos(1995), Mockus(2000), Pardalos&Romeijn(2002),
Zhigljavsky&Zilinskas(2008), Zhigljavsky(1991)}) can often deal
with the stated problems in a simpler manner than the
deterministic algorithms (being also suitable for the problems
where the evaluations of the functions are corrupted by noise).
However, there can be difficulties with these methods, as well
(e.g., in studying their convergence properties). Several restarts
can also be involved, requiring even more functions evaluations.
Moreover, solutions found by many stochastic algorithms
(especially, by popular heuristic methods like evolutionary
algorithms, simulated annealing, etc.; see, e.g.,
\cite{Michalewicz(1996), Pardalos&Romeijn(2002), Holland(1975),
Schneider&Kirkpatrick(2006), Yang(2010), Hansen&Ostermeier(2001),
Karaboga&Akay(2009), Vaz&Vicente(2007)}) can be only local
solutions to the problems, far from the global ones. This can
preclude such methods from their usage in practice, especially
when an accurate estimate of the global solution is required. That
is why we concentrate, hereafter, on deterministic approaches.

The possibility to outperform the `brute-force computation'
techniques in solving global optimization problems is
fundamentally based on the availability of some realistic a priori
assumptions characterizing the objective function and eventual
constraints (see, e.g., \cite{Horst&Pardalos(1995),
Pardalos&Romeijn(2002), Pinter(1996), Sergeyev&Kvasov(2008),
Strongin&Sergeyev(2000), Zhigljavsky&Zilinskas(2008)}). They serve
as mathematical tools for obtaining estimates of the global
solution related to a finite number of function trials and,
therefore, play a crucial role in the construction of any
efficient global search algorithm. As observed, e.g., in
\cite{Horst&Pardalos(1995), Stephens&Baritompa(1998)}, if no
particular assumptions are made on the objective function and
constraints, any finite number of function evaluations cannot
guarantee getting close to the global minimum value, since this
function may have very high and narrow peaks.

One of the natural and powerful (from both the theoretical and the
applied points of view) assumptions on the global optimization
problem is that the objective function (and constraints) has (have)
bounded slopes. In other words, any limited change in the object
parameters yields limited changes in the characteristics of the
objective performance. This assumption can be justified by the
fact that in technical systems the energy of change is always
limited (see the related discussion
in~\cite{Strongin&Sergeyev(2000)}). One of the most popular
mathematical formulations of this property is the Lipschitz
continuity condition, which assumes that the difference (in the
sense of a chosen norm) of any two function values is majorized by
the difference of the corresponding function arguments, multiplied
by a positive factor $L < \infty$. In this case, the function is
said to be Lipschitz and the corresponding factor $L$ is said to
be the Lipschitz constant. The problem involving Lipschitz
functions (the objective function and constraints) is said to be
the Lipschitz global optimization problem (see,
e.g.,~\cite{Horst&Pardalos(1995), Pinter(1996),
Sergeyev&Kvasov(2008), Sergeyev&Kvasov(2011),
Strongin&Sergeyev(2000), Zhigljavsky&Zilinskas(2008),
Kvasov&Sergeyev(2012b), Evtushenko(1985)} and the references given
therein).

The Lipschitz continuity assumption, being quite realistic for
many practical black-box problems, is also an effective tool for
obtaining accurate global optimum estimates after performing a
limited number of functions evaluations. It is used by the authors
to develop efficient and reliable deterministic methods for
solving multidimensional constrained global optimization problems
from different real-life applied areas (as, e.g., the
problem~\eqref{eq_app1}--\eqref{eq_app2}), which are characterized
by black-box multiextremal and hard to evaluate functions. In the
next Section, the Lipschitz global optimization is examined more
in detail.

\section{Lipschitz global optimization problem}

A general Lipschitz global optimization problem can be formalized
as follows (see, e.g., \cite{Pinter(1996), Sergeyev&Kvasov(2008),
Strongin&Sergeyev(2000), Zhigljavsky&Zilinskas(2008),
Evtushenko(1985)}):
 \begin{equation}
    \label{LGOP_f}
    f^* = f(x^*) = \min f(x), \hspace{5mm} x \in \Omega \subset \mathbb{R}^N,
 \end{equation}
where $\Omega$ is a bounded set defined as
 \begin{equation}
    \label{LGOP_gi}
    \Omega = \{x \in D : \, \zeta_i(x) \leq 0, \, 1 \leq i \leq p\},
 \end{equation}
 \begin{equation}
    \label{LGOP_D}
    D =[a,b\,]=\{x \in \mathbb{R} ^N : a(j)\leq x(j)\leq b(j), \, 1\leq
    j \leq N \}, \hspace{3mm} a,b \in \mathbb{R}^N,
 \end{equation}
with $N$ being the problem dimension. In
\eqref{LGOP_f}--\eqref{LGOP_D}, the objective function~$f(x)$ and
the constraints $\zeta_i(x)$, \mbox{$1 \leq i \leq p$}, are
multiextremal, non necessarily differentiable, black-box and hard
to evaluate functions that satisfy the Lipschitz condition over
the search hyperinterval $D$:
 \begin{equation}
  \label{LGOP_L}
  | f(x')-f(x'')| \leq L \|x'-x''\|, \hspace{5mm} x',x'' \in D,
 \end{equation}
 \begin{equation}
  \label{LGOP_Li}
  | \zeta_i(x')-\zeta_i(x'')| \leq L_i \|x'-x''\|,
  \hspace{5mm} x',x'' \in D, \hspace{5mm} 1 \leq i \leq p,
 \end{equation}
where $\|\cdot\|$ denotes, usually, the Euclidean norm, $L$ and
$L_i$, $1\leq i \leq p$, are the (unknown) Lipschitz constants
such that $0 < L < \infty$, $0 < L_i < \infty$, $1\leq i \leq p$.
If $p = 0$ in~(\ref{LGOP_gi}), the problem is said to be
box-constrained.

The admissible region~$\Omega$ can consist of disjoint, non-convex
subregions because of the multiextremality of the constraints
$\zeta_i(x)$. Moreover, these constraints can be partially
defined, i.e., a constraint $\zeta_{i+1}(x)$ (or the objective
function $f(x)$) can be defined only over subregions where
$\zeta_i(x) \leq 0$, $1 \leq i \leq p$ (see, e.g.,
\cite{Sergeyev&Kvasov(2011), Strongin&Sergeyev(2000)} for more
details and applied examples).

Problem~\eqref{LGOP_f}, \eqref{LGOP_D}, \eqref{LGOP_L} with a
differentiable objective function having the Lipschitz (with an
unknown Lipschitz constant) gradient $f'(x)$ (which could be
itself a multiextremal black-box function) is sometimes included
in the same class of Lipschitz global optimization problems (see,
e.g., the references given in \cite{Strongin&Sergeyev(2000),
Kvasov&Sergeyev(2012a), Kvasov&Sergeyev(2009)}).

As evidenced, e.g., in \cite{Sergeyev&Kvasov(2011),
Strongin&Sergeyev(2000)}, it is not easy to manage multiextremal
constraints (\ref{LGOP_gi}) within the context of Lipschitz global
optimization. For example, the traditional penalty approach (see,
e.g., the references in \cite{Floudas&Pardalos(2009),
Horst&Pardalos(1995), Liuzzi:et:al.(2010)}) can lead to extremely
high Lipschitz constants, thus forcing degeneration of the
methods. In this connection, a promising approach called the index
scheme (see, e.g., \cite{Strongin&Sergeyev(2000),
Barkalov&Strongin(2002), Sergeyev:et:al.(2001),
Sergeyev:et:al.(2007b)}) can be applied. It does not introduce
additional variables and/or parameters by opposition as, e.g.,
many traditional penalty approaches do, and reduces the general
constrained problem \eqref{LGOP_f}--\eqref{LGOP_Li} to a
box-constrained discontinuous one.

Therefore, in order to give an insight into the principal ideas of
the authors' techniques for solving the stated problem, box-con\-strained Lipschitz global optimization problem
(\ref{LGOP_f}), (\ref{LGOP_D}), (\ref{LGOP_L}) will be considered
in the following.

\begin{figure}
\begin{center}
\includegraphics[width=130mm,keepaspectratio]{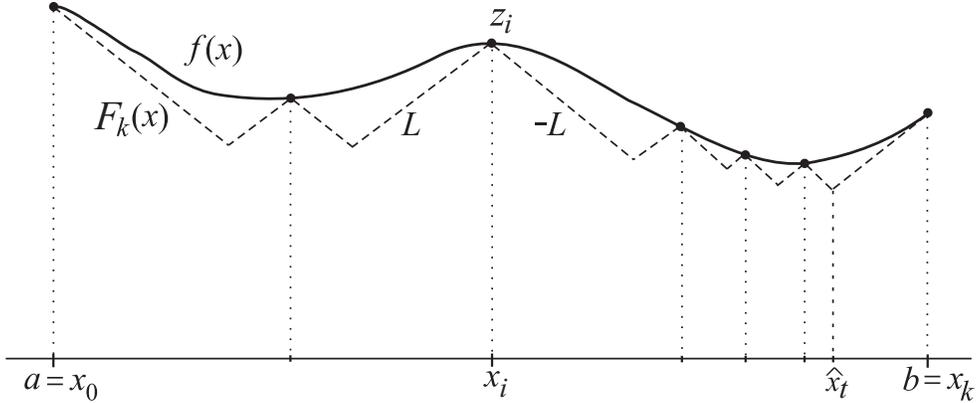}
\caption{Lower bounding function $F_k(x)$ (dashed line)
constructed for a Lipschitz function~$f(x)$ (solid line)
over~$[a,b]$ after having performed $k+1$ function trials (in this
case, $k=6$).} \label{fig:Piyavskij_L}
\end{center}
\end{figure}

Once a valid estimate of the Lipschitz constant is known and some
function trials are performed, the Lipschitz
condition~(\ref{LGOP_L}) allows us to easily find the lower bounds
of a Lipschitz function at different subregions of the search domain
$D$ from~(\ref{LGOP_D}). Let us consider, for the sake of example, a
one-dimensional objective function $f(x)$ defined over an interval
$[a,b]$ (see Figure~\ref{fig:Piyavskij_L}) that satisfies the
Lipschitz condition~(\ref{LGOP_L}) with a known Lipschitz constant
$L$. If the function values $z_i$ have been obtained at points
$x_i$, $0 \leq i \leq k$ (see black dots on the objective function
graph in Figure~\ref{fig:Piyavskij_L}), the following inequality is
satisfied over $[a,b]$:
 \begin{equation} \label{fF}
  f(x) \geq F_k(x) = \max_{0 \leq i \leq k} \{z_{i} - L|x-x_{i}|\},
 \end{equation}
where $F_k(x)$ is a piecewise linear function (called lower
bounding or minorant function, see, e.g.,
\cite{Sergeyev&Kvasov(2011), Strongin&Sergeyev(2000),
Kvasov&Sergeyev(2012b)}; its graph is drawn by a dashed line in
Figure~\ref{fig:Piyavskij_L}).

A method (e.g., the Piyavskij--Shubert method being one of the
first methods in Lipschitz global optimization, see
\cite{Horst&Pardalos(1995), Sergeyev&Kvasov(2011),
Strongin&Sergeyev(2000), Kvasov&Sergeyev(2012b),
Evtushenko(1971)}), using in its work this simple but efficient
geometric interpretation, iteratively constructs an auxiliary
function which bounds the objective function $f(x)$ from below and
evaluates $f(x)$ at a point ($\hat{x}_t$ in
Figure~\ref{fig:Piyavskij_L}) corresponding to a minimum of the
bounding function. This point is easy to find (see, e.g.,
\cite{Horst&Pardalos(1995), Sergeyev&Kvasov(2011),
Strongin&Sergeyev(2000), Kvasov&Sergeyev(2012b)}). The methods of
this type form the class of geometric algorithms that are based on
constructing, updating, and improving auxiliary piecewise
functions built by using an estimate of the Lipschitz
constant~$L$. It should be noted in this connection, that similar
ideas are used in many other surrogate-based optimization methods
(see, e.g., \cite{Shan&Wang(2010), Booker:et:al.(1999),
Jones(2001b), Forrester:et:al.(2008)}). As shown, e.g.,
in~\cite{Sergeyev&Kvasov(2008), Strongin&Sergeyev(2000)}, there
exists a strong relationship between the geometric approach and
another possible technique for solving the stated
problem\hspace{1pt}---\hspace{1pt}the so-called
information-statistical approach (see, e.g.,
\cite{Strongin&Sergeyev(2000), Strongin(1978)} and
also~\cite{Mockus(2000), Zhigljavsky&Zilinskas(2008)} for other
probabilistic techniques). Together with the geometric ideas of
the Piyavskij-Shubert method, it has consolidated foundations of
the Lipschitz global optimization.

In order to develop Lipschitz global optimization methods, the
Lipschitz constant $L$ from~(\ref{LGOP_L}) should be estimated. It
can be done in several ways. For example, the Lipschitz constant
can be given a priori (see, e.g.,~\cite{Horst&Pardalos(1995),
Evtushenko(1971), Evtushenko&Posypkin(2011)}). More practical
approaches are based on an adaptive estimation of $L$ in the
course of the search: such algorithms can use either an adaptive
global estimate of the Lipschitz constant (see,
e.g.,~\cite{Pinter(1996), Strongin&Sergeyev(2000), Strongin(1978),
Kvasov&Sergeyev(2003)}) valid for the whole search domain $D$, or
adaptive local estimates~$L_i$ valid only for some subregions $D_i
\subset D$ (see, e.g.,~\cite{Sergeyev&Kvasov(2008),
Strongin&Sergeyev(2000), Sergeyev(1995a), Sergeyev(1998a),
Kvasov:et:al.(2003)}). Estimating local Lipschitz constants during
the work of a global optimization algorithm allows one to
significantly accelerate the global search. Balancing between
local and global information must be performed in an appropriate
way (see, e.g.,~\cite{Sergeyev&Kvasov(2008),
Strongin&Sergeyev(2000), Sergeyev(1995a)}) since an unjustified
usage of local information can lead to the loss of the global
solution (see, e.g.,~\cite{Stephens&Baritompa(1998)}). Finally,
multiple estimates of $L$ can be also used (see,
e.g.,~\cite{Sergeyev&Kvasov(2008), Kelley(1999),
Jones:et:al.(1993), Sergeyev&Kvasov(2006),
Gablonsky&Kelley(2001)}). We would like to emphasize here that
either the Lipschitz constant is given and an algorithm is
developed correspondingly, or it is not known but there exist a
sufficiently large number of parameters of the considered
algorithm ensuring its convergence (convergence properties of the
Lipschitz global optimization methods are thoroughly examined,
e.g., in~\cite{Sergeyev&Kvasov(2008), Strongin&Sergeyev(2000),
Strongin(1978), Sergeyev(1998b)}).

Considering both the theoretical generality and the application
diffusion of the Lipschitz global optimization
problem~\eqref{LGOP_f}, \eqref{LGOP_D}, \eqref{LGOP_L}, it is used
by the authors to mathematically model various real-life optimal
design problems (see \cite{Sergeyev&Kvasov(2008),
Strongin&Sergeyev(2000), Golubev:et:al.(2010),
Kvasov:et:al.(2008), Kvasov:et:al.(2012), Sergeyev:et:al.(1999)}).

\section{Some deterministic tools in Lipschitz global optimization}

In this Section, some innovative deterministic approaches
developed by the authors for constructing efficient global
optimization techniques are briefly presented as in
\cite{Kvasov&Sergeyev(2012c)}. The consolidated success of these
ideas, confirmed by important international publications and
presentations around the world, allows the authors' group, on the
one hand, to develop promising optimization approaches over a solid
scientific basis, thus eliminating the theoretical faults risks,
and, on the other hand, to tackle difficult black-box practical
optimization problems (e.g., from control theory, environmental
sciences and geological mechanics, electrical engineering and
telecommunications, gravitational physics, etc.) with more
efficiency with respect to the techniques traditionally used by
the optimal design engineers.

\subsection{`Divide-the-Best' algorithms}

Many global optimization algorithms (of both deterministic and
stochastic types) have a similar structure. Therefore, several
attempts aiming to construct a general framework for describing
computational schemes and providing their convergence conditions
in a unified manner have been made (see, e.g.,
\cite{Floudas&Pardalos(2009), Horst&Pardalos(1995), Pinter(1996),
Grishagin:et:al.(1997)}). One of the more flexible and robust
among such unifying schemes is the `Divide-the-Best' approach (see
\cite{Sergeyev&Kvasov(2008), Sergeyev(1998b)}), which generalizes
both the schemes of adaptive partition~\cite{Pinter(1996)} and
characteristic~\cite{Sergeyev&Kvasov(2008),
Strongin&Sergeyev(2000), Grishagin:et:al.(1997)} algorithms,
widely used for describing and studying numerical global
optimization methods.

\begin{figure}
\begin{center}
\includegraphics[width=110mm,keepaspectratio]{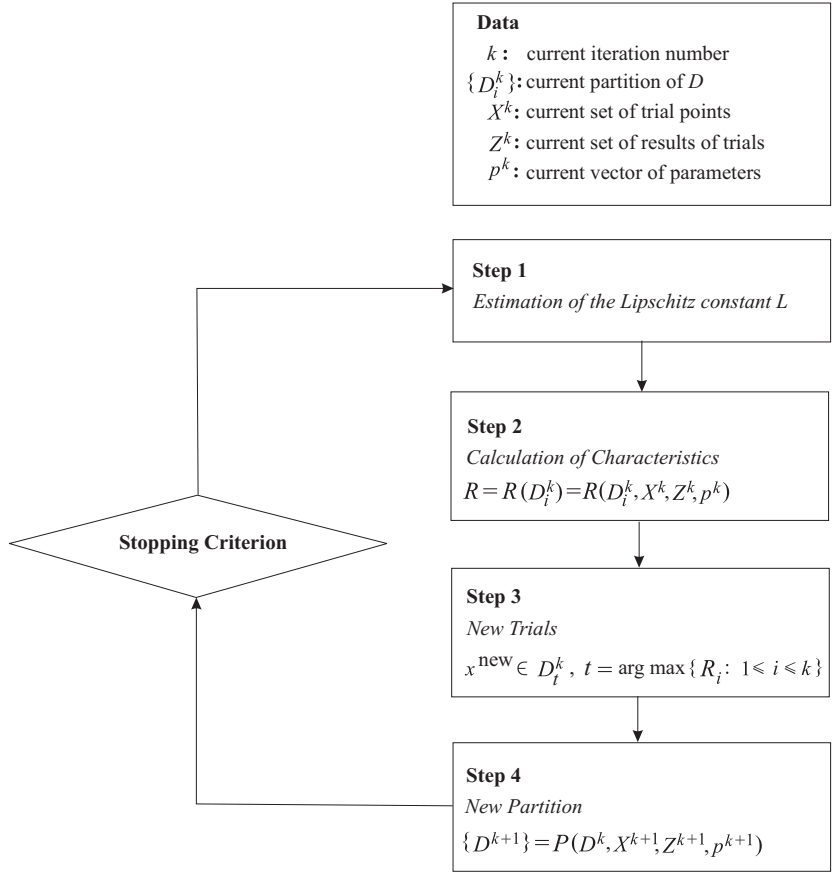}
\caption{Flow chart of `Divide-the-Best' scheme.} \label{fig:DBA}
\end{center}
\end{figure}

In this scheme (the flow chart of its generic iteration is
reported in Figure~\ref{fig:DBA}), given a vector $p$ of the
method parameters, an adaptive partition of the admissible region
$D$ from~\eqref{LGOP_D} into a collection $\{D^k_i\}$ of the
finite number of robust subsets $D^k_i$ is considered at each
iteration $k$. The `merit' (called characteristic) $R_i$ of each
subset (see Step~2 in Figure~\ref{fig:DBA}) for performing a
subsequent, more detailed, investigation (see Steps~3 and 4 in
Figure~\ref{fig:DBA}) is estimated on the basis of the obtained
information $X^k$, $Z^k$ about the objective function. The best
(in some predefined sense) characteristic obtained over some
hyperinterval $D^k_t$ corresponds to a higher possibility to find
the global minimizer within $D^k_t$ (see Step 3). This
hyperinterval is subdivided at the next iteration of the
algorithm. Naturally, more than one `promising' hyperinterval can
be partitioned at every iteration.

Several strategies (mainly, in the context of the geometric
approach) for selection of a subset for further partitioning (see
Step~3 in Figure~\ref{fig:DBA}) and for performing this
partitioning (by means of an operator $P$, see Step~4 in
Figure~\ref{fig:DBA} and the next subsection 4.2) are proposed by
the authors from a general viewpoint and successfully used for
solving practical applications (see, e.g., the references
in~\cite{Sergeyev&Kvasov(2008), Sergeyev&Kvasov(2011)}).

Regarding the stopping criteria, one can constantly check, e.g.,
the volume of a hyperinterval with the best characteristic or
depletion of computing resources such as the maximum number of
trials. The verification of a stopping criterion can be performed
at any step of the current iteration of the algorithm.

Convergence properties of the `Divide-the-Best' family for
different types of characteristic values and partition operators
are studied in~\cite{Sergeyev&Kvasov(2008), Sergeyev(1998b)}.
Great attention is given to situations (very important in
practice) when conditions of global (local) convergence are
satisfied not in the whole search domain~$D$, but only in its
small subregion (or a set of subregions). This can correspond, for
example, to Lipschitz global optimization algorithms that work
underestimating the Lipschitz constant or which are oriented on
using local information in subregions of $D$ (see, e.g.,
\cite{Sergeyev&Kvasov(2008), Strongin&Sergeyev(2000),
Sergeyev(1998b)}). It should be also noted that the described
scheme can be successfully applied to constructing parallel
multidimensional global optimization algorithms
\cite{Strongin&Sergeyev(2000), Grishagin:et:al.(1997)}.

\subsection{Efficient partitioning strategy}

Regarding the partitioning strategies (partitioning operator $P$
on Step~4 in Figure~\ref{fig:DBA}), the main attention of the
authors is focused on the diagonal partition strategies (see the
references in \cite{Pinter(1996), Sergeyev&Kvasov(2008),
Sergeyev&Kvasov(2011), Sergeyev&Kvasov(2006)}).

\begin{figure}
\begin{center}
\includegraphics[width=80mm,keepaspectratio]{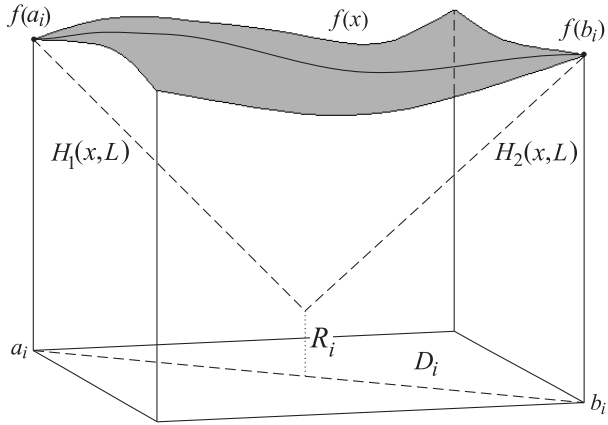}
\caption{Lower bounding a Lipschitz function $f(x)$ over a
hyperinterval $D_i$ in a diagonal algorithm.}
\label{fig:DiagonalLines}
\end{center}
\end{figure}

In this approach, the initial hyperinterval $D$ from
\eqref{LGOP_D} is partitioned into a set of smaller
hyperintervals, the objective function is evaluated only at two
vertices corresponding to the main diagonal of hyperintervals of
the current partition of $D$ (see, e.g., points $a_i$ and $b_i$ of
a hyperinterval $D_i$ in Figure~\ref{fig:DiagonalLines}), and the
results of these evaluations are used to select a hyperinterval
for the further subdivision. The diagonal approach has a number of
attractive theoretical properties and has proved to be efficient
in solving applied problems.

First, it allows one to easily perform an extension of efficient
univariate global optimization algorithms to the multidimensional
case (see, e.g., \cite{Sergeyev&Kvasov(2008),
Sergeyev&Kvasov(2011), Sergeyev&Kvasov(2006)}). In fact, in order
to calculate the characteristic $R_i$ of a multidimensional
subregion $D_i$, some one-dimensional characteristics can be used
as prototypes. After an appropriate transformation they can be
applied to the one-dimensional segment being the main diagonal
$[a_i, b_i]$ of the hyperinterval $D_i$ (see Lipschitz-based lower
bounding functions $H_1$ and $H_2$ in
Figure~\ref{fig:DiagonalLines}).

Second, the diagonal approach is close from the computational point
of view to one of the simplest strategies---centre-sampling
technique (see, e.g.,~\cite{Conn:et:al.(2009), Jones:et:al.(1993),
Evtushenko(1985), Evtushenko&Posypkin(2011),
Gablonsky&Kelley(2001)})---but at the same time, the objective
function is evaluated at two points of each subregion, providing in
this way more information about the function over the subregion than
centre-sampling methods.

Different exploration techniques based on various diagonal
adaptive partition strategies are analyzed, e.g.,
in~\cite{Sergeyev&Kvasov(2008), Sergeyev&Kvasov(2011),
Sergeyev(2000)}. It is demonstrated that partition strategies
traditionally used in the framework of the diagonal approach do
not fulfil the requirements of computational efficiency because of
the execution of many redundant trials. Such a redundancy slows
down significantly the global search in the case of costly
functions.

An efficient diagonal partition strategy is therefore proposed
in~\cite{Sergeyev&Kvasov(2008), Sergeyev(2000)}, that allows one
to avoid the computational redundancy of traditional diagonal
schemes. In contrast to these schemes, the proposed strategy produces
regular meshes of the function evaluation points in such a way
that one vertex where~$f(x)$ is evaluated can belong to several
hyperintervals (up to $2^N$, $N$~is the problem dimension
from~(\ref{LGOP_D})). Thus, the time-consuming procedure of the
function evaluations is replaced by a significantly faster
operation of reading (up to $2^N$ times) the function values
obtained at the previous iterations and saved in a special
database (see, e.g., \cite{Kvasov(2008a), Kvasov(2008b)}). Hence,
this partition strategy considerably speeds up the search and
also leads to saving computer memory. It is particularly important
that these advantages become more pronounced when the problem dimension $N$ increases (see, e.g., \cite{Sergeyev&Kvasov(2008), Sergeyev&Kvasov(2006),
Kvasov&Sergeyev(2003)}).

A novel scheme for creating fast Lipschitz global optimization
algorithms is, thus, introduced by the authors. It relies on the
efficient diagonal partition strategy allowing an efficient
extension of popular one-dimensional Lipschitz global optimization
algorithms to the multidimensional case. In a sense, this scheme
combines the ideas of the diagonal approach and Peano
space-filling curves (see, e.g.,~\cite{Strongin&Sergeyev(2000),
Strongin(1978), Sergeyev:et:al.(2013)}). Innovative
multidimensional diagonal algorithms for solving Lipschitz global
optimization problems, based on different ways for obtaining the
Lipschitz information and developed in the framework of the
efficient diagonal scheme, are proposed by the authors and their
convergence properties are analyzed, e.g.,
in~\cite{Sergeyev&Kvasov(2008), Sergeyev&Kvasov(2006),
Kvasov&Sergeyev(2003)}.

\subsection{Balancing local and global information}

Is well known (see, e.g., \cite{Floudas&Pardalos(2009),
Horst&Pardalos(1995), Sergeyev&Kvasov(2008),
Strongin&Sergeyev(2000), Stephens&Baritompa(1998)}) that the usage
of the only global information on the objective function and
constraints during optimization can lead to a slow convergence of
algorithms to global minimizers. Therefore, particular attention
is paid by the authors to the usage of local information in
global optimization methods, as well. One of the traditional ways
in this context (see, e.g., \cite{Floudas&Pardalos(2009),
Horst&Pardalos(1995), Pardalos&Romeijn(2002)}) recommends stopping
the global procedure and switching to a local optimization method
in order to improve the solution and to accelerate the search
during its final phase. Unfortunately, applying this technique can
lead to some problems related to the combination of global and
local phases, the main problem being that of determining when to
stop the global procedure and start the local one. A premature
arrest can provoke the loss of the global solution whereas a late
one can slow down the search.

Theoretical and experimental results obtained by the authors (see,
e.g., \cite{Sergeyev&Kvasov(2008), Strongin&Sergeyev(2000),
Sergeyev(1995a), Kvasov:et:al.(2003)}) confirm that more fruitful
approaches can be considered. The first one is the so-called local
tuning approach~\cite{Sergeyev(1995a)} allowing global
optimization algorithms to tune their behaviour to the shape of
the functions at different parts of the search domain by
estimating the local Lipschitz constants.

In fact, the Lipschitz constant $L$ has a significant influence on
the convergence speed of the Lipschitz global optimization
algorithms and the problem of its specifying is of great
importance. Accepting, for instance, too high a value of $L$ for a
concrete objective function means assuming that the function has
complicated structure with sharp peaks and narrow attraction
regions of minimizers within the whole admissible region. Thus, if
the value of $L$ does not correspond to the real behaviour of the
objective function, it can lead to a slow convergence of the
algorithm to the global minimizer. Global optimization algorithms
using in their work a global estimate of $L$ (or some values of
$L$ given a priori) do not take into account local information
about behaviour of the objective function over every small
subregion of $D$. Therefore, estimating local Lipschitz constants
allows one to significantly accelerate the global search (see,
e.g., \cite{Sergeyev&Kvasov(2008), Strongin&Sergeyev(2000),
Sergeyev(1998a), Kvasov:et:al.(2003)}).

The second technique regards a continual local improvement of the
current best solution incorporated in a global search procedure
(see, e.g.,~\cite{Sergeyev&Kvasov(2008), Sergeyev&Kvasov(2006),
Lera&Sergeyev(2010b), Lera&Sergeyev(2013)}). Particularly, it
forces the global optimization method to make a local improvement
of the best approximation of the global minimum immediately after
a new approximation better than the current one is found. These
techniques become even more efficient when information about the
objective function derivatives is available (see,
e.g.,~\cite{Kvasov&Sergeyev(2012a), Kvasov&Sergeyev(2009)}).

\subsection{Computational aspects}

A particular attention is paid by the authors to the problem of
testing global optimization algorithms. As widely accepted, a set
of test functions is usually taken for this purpose, problems from
this set are solved by the algorithms to be compared, and a
conclusion about the efficiency of these algorithms is made on the
basis of the obtained numerical results. This approach, being an
important instrument for acquiring a knowledge about the existing
and new global optimization algorithms, presents at the same time
some limitations since the conclusions made can be valid only for
the selected functions, and their propagation to a more wide set
of functions requires particular caution. Testing an algorithm on
a relatively large set of test functions can, in a sense, diminish
these limitations, but it needs, among other things, the coding of
the functions, and it is a tedious and time-consuming job.
Moreover, the lack of such information as number of local optima,
their locations, attraction regions, local and global values,
describing global optimization tests taken from real-life
applications, creates additional difficulties in verifying
validity of the algorithms. Therefore, the global optimizers are
very interested in simple and powerful software tools realizing
test problems. As observed, e.g., in~\cite{Sergeyev&Kvasov(2008),
Strongin&Sergeyev(2000), Zhigljavsky&Zilinskas(2008),
Rios&Sahinidis(2013), Grishagin(1978), Dolan&More(2002),
More&Wild(2009)}, a well designed testing framework is of the
primary importance in identifying the merits of each algorithm and
implementation.

To tackle the problem of testing global optimization algorithms
systematically, the GKLS-generator described in
\cite{Gaviano:et:al.(2003)} is proposed by the authors' group. The
generator produces several classes of multidimensional and
multiextremal test functions with known local and global minima.
Each test class provided by the generator includes 100 functions. By
changing the user-defined parameters, classes with different
properties can be created. For example, fixed dimension of the
functions and number of local minima, a more difficult class can be
created either by shrinking the attraction region of the global
minimizer, or by moving the global minimizer closer to the domain
boundary.

The generator is available on the ACM Collected Algorithms (CALGO)
database (the CALGO is part of a family of publications produced
by the Association for Computing Machinery) and it is also
downloadable for free from
\url{http:\\wwwinfo.dimes.unical.it\~yaro\GKLS.html}. It has
already been downloaded by companies and research organizations
from more than 40 countries of the world.

\section{Some numerical results}

To conclude, we would like to report some numerical results
obtained by using a Lipschitz global optimization method
proposed by the authors in~\cite{Sergeyev&Kvasov(2006)}. In
developing this method for solving problem~\eqref{LGOP_f},
\eqref{LGOP_D}, \eqref{LGOP_L}, techniques from the previous
Section have been applied. Particularly, it is a multidimensional
`Divide-the-Best' global optimization method that uses in its work
multiple estimates of the Lipschitz constant and based on
efficient diagonal partitions.

Numerical results performed on the GKLS-generator to compare this algorithm with two algorithms belonging to the same class of
methods for solving problem~\eqref{LGOP_f}, \eqref{LGOP_D},
\eqref{LGOP_L}
--- the DIRECT algorithm from~\cite{Jones:et:al.(1993)} and its
locally-biased modification DIRECT{\it l} \hspace{1mm}from
\cite{Gablonsky&Kelley(2001)} --- are presented here, as described
in~\cite{Sergeyev&Kvasov(2006)}. As known, both of these methods
are widely used in solving practical engineering problems (see,
e.g., the references in~\cite{Floudas&Pardalos(2009),
Kelley(1999), Sergeyev&Kvasov(2006)}). Moreover, as shown
numerically in~\cite{Grbic:et:al.(2013)}, they often outperform
metaheuristic algorithms as, e.g., the Firefly algorithm
(see~\cite{Yang(2010)}) belonging to the widely used family of
Particle Swarm Optimization algorithms (see,
e.g.,~\cite{Yang(2010), Vaz&Vicente(2007)}).

Eight GKLS classes of continuously differentiable test functions
of dimensions $N=2$, 3, 4, and 5 have been used. For each
dimension, both a `hard' and a `simple' classes have been
considered. The difficulty of a class was increased either by
decreasing the radius of the attraction region of the global
minimizer, or by decreasing the distance from the global minimizer
$x^*$ to the domain boundaries.

The global minimizer $x^* \in D$ was considered to be found when
the algorithm generated a trial point~$x'$ inside a hypercube with
a vertex $x^*$ and the volume smaller than the volume of the
initial hypercube $D=[a,b]$ multiplied by an accuracy
coefficient~$\Delta$, $0< \Delta \leq 1$, i.e.,
 \begin{equation} \label{Delta}
  |x'(j) - x^*(j) | \leq \sqrt[N]{\Delta}(b(j)-a(j))
 \end{equation}
for all $i$, $1 \leq j \leq N$, where $N$ is from~\eqref{LGOP_D}.
The algorithm stopped either when the maximal number of trials
equal to 1 000 000 was reached, or when condition~\eqref{Delta}
was satisfied.

In view of the high computational complexity of each trial of the
objective function, the methods were compared in terms of the
number of evaluations of $f(x)$ required to satisfy
condition~\eqref{Delta}. The number of hyperintervals generated
until condition~\eqref{Delta} is satisfied, was taken as the
second criterion for comparison of the methods. This number
reflects indirectly degree of qualitative examination of $D$
during the search for a global minimum (see, e.g.,
\cite{Sergeyev&Kvasov(2008), Sergeyev&Kvasov(2006),
Kvasov&Sergeyev(2003)}).

Results of numerical experiments with eight GKLS tests classes are
reported in Tables \ref{table1}--\ref{table2}. These tables show,
respectively, the maximal number of trials and the corresponding
number of generated hyperintervals required for satisfying
condition~\eqref{Delta} for a half of the functions of a
particular class (columns ``50\%'') and for all 100 function of
the class (columns ``100\%''). The notation ``$>$ 1 000 000
$(k)$'' means that after 1 000 000 trials the method under
consideration was not able to solve $k$ problems.

Note that on a half of test functions from each class (which were
simple for each method with respect to the other functions of the
class) the algorithm from~\cite{Sergeyev&Kvasov(2006)} manifested
a good performance with respect to DIRECT and DIRECT{\it l} in
terms of the number of generated trial points (see
Table~\ref{table1}). When all functions were taken in
consideration, the number of trials produced by the new algorithm
was significantly fewer in comparison with two other methods (see
columns ``100\%'' of Table~\ref{table1}), providing at the same
time a good examination of the admissible region (see
Table~\ref{table2}).

\begin{table}[t]
\begin{center}
\caption{Number of trial points for 800 GKLS test functions.}
\label{table1} \scriptsize
\begin{tabular}{|c|c|c|r|r|r|r|r|r|}\hline
$N$ & $\Delta$ & Class & \multicolumn{3}{|c|}{50\%} &
\multicolumn{3}{|c|}{100\%}\\
\cline{4-9}& & &DIRECT &DIRECT{\it l} &New &DIRECT &DIRECT{\it l} &New\\
\hline
2 &$10^{-4}$ & simple &111  &152 &166  &1159  &2318 &403\\
2 &$10^{-4}$ & hard &1062  &1328 &613  &3201  &3414 &1809\\
\hline
3 &$10^{-6}$ & simple &386   &591 &615  &12507 &13309 &2506\\
3 &$10^{-6}$ & hard &1749  &1967 &1743 &$>$1000000 (4) &29233 & 6006\\
\hline
4 &$10^{-6}$ & simple &4805  &7194 &4098  &$>$1000000 (4) &118744 &14520\\
4 &$10^{-6}$ & hard &16114 &33147 &15064 &$>$1000000 (7) &287857 &42649\\
\hline
5 &$10^{-7}$ & simple &1660  &9246 &3854  &$>$1000000 (1) &178217 &33533\\
5 &$10^{-7}$ & hard &55092  &126304 &24616  &$>$1000000 (16) &$>$1000000 (4) &93745\\
\hline
\end{tabular}
\end{center}

\begin{center} \caption{Number of hyperintervals for 800 GKLS test
functions.} \label{table2} \scriptsize
\begin{tabular}{|c|c|c|r|r|r|r|r|r|}\hline
$N$ & $\Delta$ & Class & \multicolumn{3}{|c|}{50\%} &
\multicolumn{3}{|c|}{100\%}\\
\cline{4-9}& & &DIRECT &DIRECT{\it l} &New &DIRECT &DIRECT{\it l} &New\\
\hline
2 &$10^{-4}$ & simple &111  &152 &269  &1159  &2318 & 685\\
2 &$10^{-4}$ & hard &1062 &1328 &1075 &3201 &3414 &3307 \\
\hline
3 &$10^{-6}$ & simple &386  &591  &1545 &12507 &13309 &6815\\
3 &$10^{-6}$ & hard &1749 &1967 &5005 &$>$1000000 &29233 &17555\\
\hline
4 &$10^{-6}$ & simple &4805 &7194 &15145 &$>$1000000 &118744 &73037\\
4 &$10^{-6}$ & hard &16114 &33147 &68111 &$>$1000000&287857 &211973\\
\hline
5 &$10^{-7}$ & simple &1660 &9246 &21377 &$>$1000000&178217 &206323\\
5 &$10^{-7}$ & hard &55092  &126304 &177927  &$>$1000000 &$>$1000000 &735945\\
\hline
\end{tabular}
\end{center}
\end{table}

\begin{figure}
\begin{center}
\includegraphics[width=130mm,keepaspectratio]{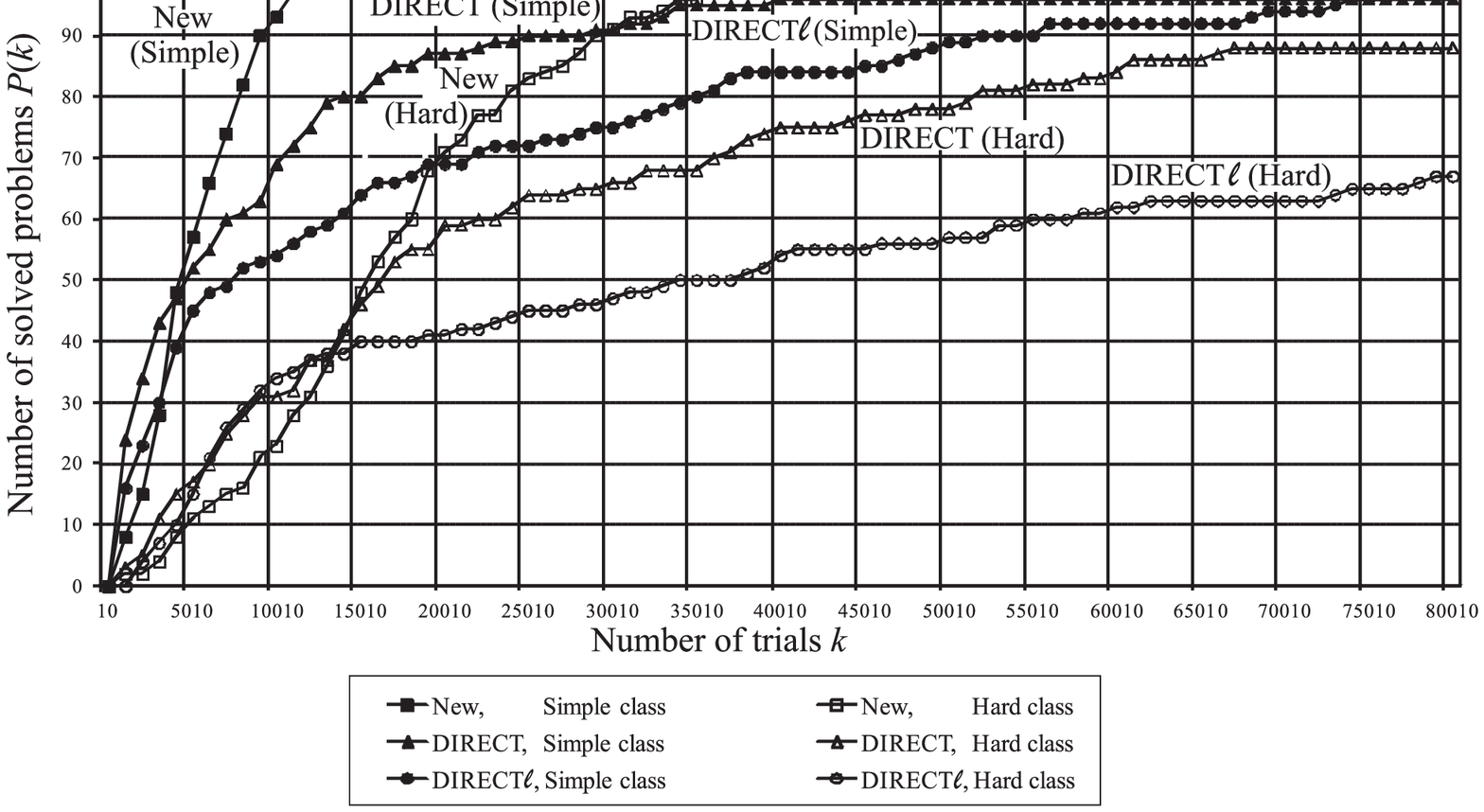}
\caption{Operating characteristics of the methods New
from~\cite{Sergeyev&Kvasov(2006)}, DIRECT, and DIRECT{\it l} on
the `simple' and `hard' four-dimensional GKLS test classes}
\label{fig:Oper}
\end{center}
\end{figure}

As it can be seen from Tables~\ref{table1}--\ref{table2}, the
method~\cite{Sergeyev&Kvasov(2006)} demonstrates a quite
satisfactory performance with respect to popular
DIRECT~\cite{Jones:et:al.(1993)} and DIRECT{\it l}
\cite{Gablonsky&Kelley(2001)} methods when multidimensional
functions with a really complex structure are minimized. Its
superiority can be also confirmed by the so-called operating
characteristics (introduced in 1978 in~\cite{Grishagin(1978)},
see~\cite{Strongin&Sergeyev(2000)} for their English language
description; they can be considered as predecessors of
`performance profiles' from~\cite{Dolan&More(2002)} and `data
profiles' from~\cite{More&Wild(2009)}). The operating
characteristics, as an indicator of the efficiency of an
optimization method, are formed by the pairs $(k, P(k))$ where $k$
($k >0$) is the number of trials and $P(k)$ ($0 \leq P(k) \leq M$)
is the number of test problems (among a set of $M$ tests) solved
by the method with less than or equal to $k$ function trials. It
is convenient to represent this indicator in a graph where each
pair (for increasing values of $k$) corresponds to a point on the
plane.

For example, Figure~\ref{fig:Oper} illustrates the operating
characteristics for the new method and the DIRECT and DIRECT{\it
l} methods while solving $M=100$ four-dimensional functions of
both `simple' and `hard' GKLS classes. The transition from the
`simple' class to the `hard' one can be clearly traced in the
diagram: e.g., the new method has solved 50 problems from the
`simple' class after 4098 trials while the solution of 50 problems
from the `hard' class has required 15064 trials (see the
intersection of the horizontal line $P(k) = 50$ with graphs
labelled `New (Simple)' and `New (Hard)' in
Figure~\ref{fig:Oper}).

By examining the operating characteristics in
Figure~\ref{fig:Oper}, it can be seen that for functions with
simple structure (up to 50 functions of the `simple' class and up
to 40 functions of the `hard' class) all three methods behave
similarly. But from the global optimization viewpoint such
functions are not interesting because it is possible to
successfully minimize them even by very naive methods. The
situation is changed when problems with complex structure should
be solved: here, both the DIRECT and DIRECT{\it l} methods
experience serious difficulties. For example, on the `hard' class
the new method has solved all 100 problems after about 43000
trials while the DIRECT method has solved 75 problems after the
same number of trials and the DIRECT{\it l} method
--- only 55 problems. Note also that even after 80000 trials
neither DIRECT nor DIRECT{\it l} methods were able to solve all
the problems of both the `simple' and `hard' classes of dimension
$N = 4$ (see the most right vertical line in
Figure~\ref{fig:Oper}). A similar situation occurs for the
operating characteristics on classes of dimensions $N = 2$, $3$,
and $5$, thus, confirming the efficiency of the proposed method
compared with the DIRECT and DIRECT{\it l} methods when solving
multidimensional multiextremal problems.

This method, as an example of the deterministic techniques
mentioned in the previous Sections, not only has manifested a high
performance on a large set of tests, but has been also
successfully applied for solving real-world global optimization
problems. For example, its application to a control theory problem
has been considered in~\cite{Kvasov:et:al.(2008)}. This problem
regards global tuning of fuzzy power system stabilizers present in
a multi-machine power system in order to damp the power system
oscillations. Power system stabilizers with conventional industry
structure are extensively used in modern power systems as an
efficient means of damping power. Traditionally their parameters
are determined by a local tuning procedure based on a
single-machine infinite-bus system in which the effects of
inter-machine and inter-area dynamics are usually ignored.
Heuristic methods (like genetic algorithms) are usually used for
their optimizing (as in many other engineering contexts) that
often leads to very rough solutions (see, e.g., the references
in~\cite{Kvasov:et:al.(2008)}). To improve overall system dynamic
performance, novel global optimization techniques have been
therefore applied by the authors' group
in~\cite{Kvasov:et:al.(2008)}.

\begin{figure}[ht]
\begin{center}
\includegraphics[width=100mm,keepaspectratio]{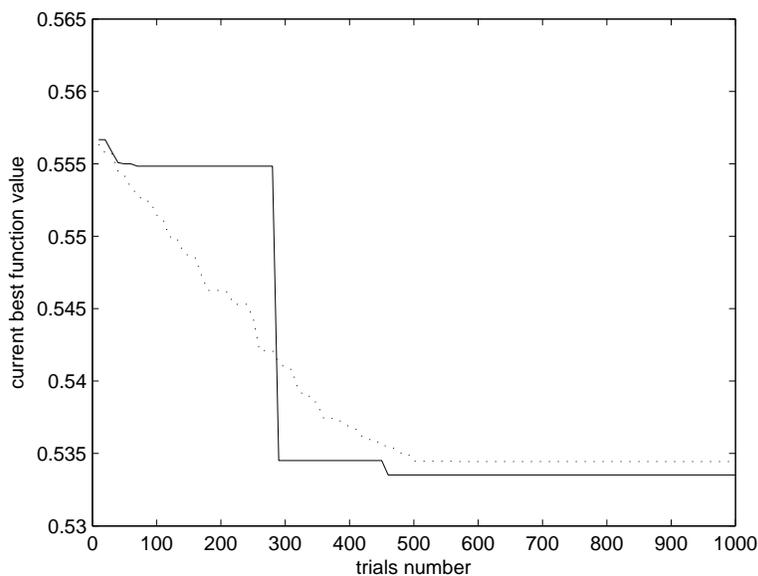}
\caption{Solutions to the problem of global tuning fuzzy power
system stabilizers (see~\cite{Kvasov:et:al.(2008)}) obtained by
applying the method~\cite{Sergeyev&Kvasov(2006)} based on the
authors' techniques (solid line) and by a traditionally used
genetic approach (dotted line).} \label{fig:GA}
\end{center}
\end{figure}

In Figure~\ref{fig:GA}, the graph that illustrates the best
solution (the axis of ordinates) obtained by a particular genetic
algorithm (often used by engineers from the control field) and by
the method~\cite{Sergeyev&Kvasov(2006)} after a number of
simulations (the axis of abscissas) is reported. It can be seen
that the global optimization method proposed by the authors spent
more function trials (namely, 284 trials) than the genetic
algorithm at the initial iterations. This phase corresponds to the
initial exploration of the search domain and it is necessary for
all global optimization techniques. On the initial phase of the
work (less than 300 trials) the genetic algorithm has found local
solutions to the problem better than those found by the
method~\cite{Sergeyev&Kvasov(2006)}, but far from the final global
solution ($f^* \approx 0.533$). However, it is more important and
should be underlined that the method~\cite{Sergeyev&Kvasov(2006)}
has determined a solution to the problem very close to the global
optimal one (as demonstrated in~\cite{Kvasov:et:al.(2008)}) in
almost half of the simulations with respect to the genetic
algorithm (284 trials for the method~\cite{Sergeyev&Kvasov(2006)}
and 500 for the genetic algorithm). Moreover, it has found an
attraction region of a new minimizer with a much better solution
to the problem (see the graph jump in Figure~\ref{fig:GA} around
450 trials) than that found by the genetic approach. Thus, when a
reasonable limit of function trials is given, the considered
method~\cite{Sergeyev&Kvasov(2006)} can determine a good estimate
of the global solution to the studied control theory problem
faster than the traditionally used genetic techniques.

Therefore, global optimization techniques briefly presented in
this survey can provide the scientists and engineers with
comprehensive and powerful tools for successful solving
challenging decision-making problems from different real-life
application areas, which are characterized by black-box
multiextremal and hard to evaluate functions. A more detailed and
systematic comparison of the described deterministic approaches
with some heuristic nature inspired techniques widely used in
engineering applications could be an interesting and useful
direction of future research.

\section*{Acknowledgements}

The authors work was supported by the project 14.B37.21.0878 of
the Ministry of Education and Science of the Russian Federation
and by the grant 1960.2012.9 awarded by the President of the
Russian Federation for supporting the leading research groups.

The authors would like to thank two anonymous referees for their
useful comments and suggestions.

%
%
%


\bibliographystyle{plain}
\bibliography{Kvasov_CST2012}

\end{document}